\documentclass{amsart} 

\usepackage{color}

\usepackage{amssymb}

\numberwithin{equation}{section} 
 
\newtheorem{thm}{Theorem}[section]

\newtheorem{prop}[thm]{Proposition}

\theoremstyle{definition}
\newtheorem{defn}{Definition}[section]
\newtheorem{exmp}{Example}[section]

\theoremstyle{remark}
\newtheorem{rem}{Remark}[section]
\newtheorem*{notation}{Notation}
\newtheorem*{Ack}{Acknowledgment}

 \DeclareMathOperator{\LL}{L}

\newcommand{\braket}[2]{\left\langle{\,{#1}\,,\,{#2}\,}\right\rangle}

\newcommand{\Z}{\mathbb{Z}}

\begin{document}

\bibliographystyle{plain}

\title{Graded Poisson Algebras}

\begin{abstract}
This note is an expanded and updated version of our entry with the same title for the 2006 Encyclopedia of Mathematical Physics. We give a brief overview of graded Poisson algebras, their main properties and their main applications, in the contexts of super differentiable and of derived algebraic geometry.
\end{abstract}

\author[A.~S.~Cattaneo]{Alberto S.~Cattaneo}
\address{Institut f\"ur Mathematik --- Universit\"at Z\"urich--Irchel ---
Winterthurerstrasse 190 --- CH-8057 Z\"urich --- Switzerland}
\email{cattaneo@math.uzh.ch}

\author[D.~Fiorenza]{Domenico~Fiorenza}
\address{Dipartimento di Matematica ``G. Castelnuovo'' --- 
Universit\`a di Roma ``La Sapienza'' ---  Piazzale Aldo Moro, 2 
--- I-00185 Roma --- Italy}
\email{fiorenza@mat.uniroma1.it}

\author[R. Longoni]{Riccardo~Longoni}  
\address{Dipartimento di Matematica ``G. Castelnuovo'' --- 
Universit\`a di Roma ``La Sapienza'' ---  Piazzale Aldo Moro, 5 --- 
I-00185 Roma --- Italy}  
\email{longoni@mat.uniroma1.it}

\thanks{A.~S.~C.~acknowledges partial support of SNF Grant No.~200020-107444/1}



\maketitle  

\tableofcontents

\section{Definitions}

\subsection{Graded vector spaces}

By a $\Z$-graded vector space (or simply, graded vector 
space) we mean a direct sum $A=\oplus_{i\in\Z}A_i$ of 
vector spaces over a field $k$ of characteristic zero. 
The $A_i$ are called the components of $A$ of degree $i$ and 
the degree of a homogeneous element $a\in A$ is denoted by $|a|$.
We also denote by $A[n]$ the graded vector space with degree 
shifted by $n$, namely, $A[n]=\oplus_{i\in\Z}(A[n])_i$ with 
$(A[n])_i=A_{i+n}$. The tensor product of two graded vector 
spaces $A$ and $B$ is again a graded vector space whose 
degree $r$ component is given by 
$(A\otimes B)_r=\oplus_{p+q=r}A_p\otimes B_q$.

The symmetric and exterior algebra of a graded vector space $A$ 
are defined respectively as $S(A)=T(A)/I_S$ and 
$\bigwedge(A)=T(A)/I_\wedge$, where $T(A)=\oplus_{n\ge 0} A^{\otimes n}$ 
is the tensor algebra of $A$ and $I_S$ (resp.\ $I_\wedge$) is the 
two-sided ideal generated by elements of the form 
$a\otimes b-(-1)^{|a|\,|b|} b\otimes a$ (resp.\ 
$a\otimes b + (-1)^{|a|\,|b|} b\otimes a$), 
with $a$ and $b$ homogeneous elements of $A$.  
The images of $A^{\otimes n}$ in $S(A)$ and $\bigwedge(A)$ 
are denoted by $S^n(A)$ and $\bigwedge^n(A)$ respectively. Notice 
that there is a canonical decalage isomorphism
$S^n(A[1])\simeq\bigwedge^n(A)[n]$.

\subsection{Graded algebras and graded Lie algebras}

We say that $A$ is a graded algebra (of degree zero) if $A$ is a 
graded vector space endowed with a degree zero bilinear associative 
product $\cdot\colon A\otimes A\to A$. A graded algebra is graded 
commutative if the product satisfies the condition 
$$a\cdot b = (-1)^{|a|\,|b|} b\cdot a$$
for any two homogeneous elements $a,b\in A$ of degree $|a|$ and $|b|$ 
respectively.

A graded Lie algebra of degree $n$ is a graded vector 
space $A$ endowed with a graded Lie bracket on $A[n]$. Such a
bracket can
 be seen as a degree $-n$ Lie bracket on $A$, i.e., as bilinear 
operation
$\{\cdot,\cdot\}\colon A\otimes A\to A[-n]$ satisfying graded
antisymmetry and graded Jacobi relations:
\begin{align*}
\{a,b\} &= -(-1)^{(|a|+n)(|b|+n)}\{b,a\}\\
\{a,\{b,c\}\} &= \{\{a,b\},c\} + (-1)^{(|a|+n)(|b|+n)}
\{b,\{a,c\}\}
\end{align*}

\subsection{Graded Poisson algebra}

We can now define the main object of interest of this note:

\begin{defn}
A graded Poisson algebra of degree $n$, or $n$-Poisson algebra, is a triple 
$(A,\cdot,\{,\})$ consisting of a graded vector space $A=\oplus_{i\in\Z}A_i$ 
endowed with a degree zero graded commutative product 
and with a degree $-n$ Lie bracket.
The bracket is required to be a biderivation of the product, namely:
\[
\{a,b\cdot c\}=\{a,b\}\cdot c + (-1)^{|b|(|a|+n)} b\cdot \{a,c\}. 
\]
\end{defn}

\begin{notation}
Graded Poisson algebras of degree zero are called Poisson algebras,
while for $n=1$ one speaks of Gerstenhaber algebras \cite{Gerst} or
of Schouten algebras.
\end{notation}

Sometimes a $\mathbb Z_2$-grading is used instead of a $\mathbb Z$-grading.
In this case, one just speaks of even and odd Poisson algebras.

\begin{exmp}
Any associative algebra can be seen as a Poisson algebra
with the trivial Lie structure, and any graded Lie algebra
can be seen as a Poisson algebra with the trivial product.
\end{exmp}

\begin{exmp}
\label{exmp:smooth}
The most classical example of a Poisson algebra is the algebra of smooth functions on ${\mathbb R}^{2n}$
endowed with usual multiplication and with the Poisson
bracket (already considered
by Poisson himself)
$\{f,g\}=\partial_{q^i}f\partial_{p_i}g-
\partial_{q^i}g\partial_{p_i}f$, where the $p_i$'s and the $q^i$'s,
for $i=1,\dots, n$, are coordinates on ${\mathbb R}^{2n}$. The
bivector field $\partial_{q^i}\wedge\partial_{p_i}$ is
induced by the symplectic form $\omega=\mathrm{d}p_i\wedge \mathrm{d}q^i$. An
immediate generalization of this example is the algebra of
smooth functions on a symplectic manifold $({\mathbb R}^{2n},\omega)$
with
the Poisson bracket
$\{f,g\}=\omega^{ij}\partial_if\partial_jg$, where
$\omega^{ij}\partial_i\wedge\partial_j$ is the bivector
field defined by the inverse of
the symplectic form
$\omega=\omega_{ij}\mathrm{d}x^i\wedge \mathrm{d}x^j$; viz.: $\omega_{ij}\omega^{jk}=\delta_i^k$.

A further generalization is when the bracket on 
$\mathcal C^\infty(\mathbb R^m)$
is defined by
$\{f,g\}=\alpha^{ij}\partial_if\partial_jg$, with the matrix function $\alpha$
not necessarily nondegenerate. The bracket is Poisson if and only if
$\alpha$ is skewsymmetric and satisfies
\[
\alpha^{ij}\partial_i\alpha^{kl}+
\alpha^{il}\partial_i\alpha^{jk}+
\alpha^{ik}\partial_i\alpha^{lj}=0.
\]

An example of this, already considered by Lie in \cite{SLie}, is 
$\alpha^{ij}(x)=f^{ij}_kx^k$, where the $f^{ij}_k$'s are the structure
constants of some Lie algebra.
\end{exmp}

\begin{exmp}\label{exmp:smoothnocoor}
Example~\ref{exmp:smooth} can be generalized to any symplectic manifold 
$(M,\omega)$. To every function $h\in\mathcal C^\infty(M)$ one associates
the Hamiltonian vector field $X_h$ which is the unique vector field satisfying
$\iota_{X_h}\omega=\mathrm{d} h$. The Poisson bracket of two functions $f$ and $g$
is then defined by
\[
\{f,g\} = \iota_{X_f}\iota_{X_g} \omega.
\]
In local coordinates, the corresponding Poisson bivector field 
is related to
the symplectic form as in Example~\ref{exmp:smooth}.

A generalization is the algebra of smooth functions on a
manifold $M$ with bracket
$\{f,g\}=\langle\alpha,\mathrm{d}f\wedge \mathrm{d}g\rangle$, where
$\alpha$ is a 
bivector field (i.e., a section of $\bigwedge^2TM$)
such that
$\{\alpha,\alpha\}_{SN}=0$, where $\{\cdot,\cdot\}_{SN}$ is the
Schouten--Nijenhuis bracket (see subsection~\ref{ssec:sn} below for
details, and
Example~\ref{exmp:smooth} for the local coordinate expression).
Such a bivector field is called a Poisson bivector field and the manifold
$M$ is called a Poisson manifold. Observe that a Poisson algebra structure on
the algebra of smooth
functions on a smooth manifold is necessarily defined this way. 
In the symplectic case,
the bivector field corresponding to
the Poisson bracket is the inverse of the symplectic form (regarded as a bundle
map $TM\to T^*M$).

The linear case described at the end of example~\ref{exmp:smooth}
corresponds to $M=\mathfrak g^*$ where $\mathfrak g$ is a (finite dimensional)
Lie algebra. The Lie bracket $\bigwedge^2 \mathfrak g\to\mathfrak g$ is
regarded as an element of 
$\mathfrak g\otimes \bigwedge^2 \mathfrak g^*\subset\Gamma(\bigwedge^2T\mathfrak g^*)$ 
and
reinterpreted as a Poisson bivector field on $\mathfrak g^*$. The Poisson algebra
structure restricted to polynomial functions is described at the beginning
of subsection~\ref{ssec:sn}.
\end{exmp}

\subsection{Batalin--Vilkovisky algebras}

When $n$ is odd, a generator for the bracket of an $n$-Poisson algebra $A$ is a
degree 
$-n$ linear map from $A$ to itself, $$\Delta\colon A\to A[-n],$$ such that  
\[
\Delta(a\cdot b)=\Delta(a)\cdot b + (-1)^{|a|}a\cdot\Delta(b) + 
(-1)^{|a|}\{a,b\}.
\]
A generator $\Delta$ is called exact if and only if it satisfies the condition 
$\Delta^2=0$, and in this case $\Delta$ becomes a derivation of the bracket:
\[
\Delta(\{a,b\})=\{\Delta(a),b\} + (-1)^{|a|+1}\{a,\Delta(b)\}.
\]

\begin{rem}
Notice that not every odd Poisson algebra $A$ admits a generator. For 
instance, a
nontrivial odd Lie algebra seen as an odd Poisson algebra with trivial
multiplication admits no generator. Moreover, even if a generator $\Delta$
for an odd Poisson algebra exists, it is far from being unique. In fact, all
different generators are obtained by adding  to
$\Delta$ a derivation of $A$ of degree $-n$. 
\end{rem}

\begin{defn}
An $n$-Poisson algebra $A$ is called an $n$-Batalin--Vilkovisky algebra,  
if it is endowed with an exact generator. 
\end{defn}

\begin{notation}
When $n=1$ it is customary to speak of Batalin--Vilkovisky algebras, or 
simply BV algebras; see \cite{B-V,Getz,Kosz}.
\end{notation}

There exists a characterization of $n$-Batalin--Vilkovisky algebras
in terms of the product and the generator only
\cite{Getz,Kosz}.  Suppose in fact that a graded vector
space $A$ is endowed with a  degree zero graded commutative
product and a  linear map
$\Delta\colon A\to A[-n]$ such that $\Delta^2=0$, satisfying the 
following ``seven-term'' relation:
\begin{multline*}
\Delta(a\cdot b\cdot c) + \Delta(a)\cdot b\cdot c + (-1)^{|a|} a\cdot
\Delta(b)\cdot c + (-1)^{|a|+|b|} a\cdot b\cdot \Delta(c) =\\ 
= \Delta(a\cdot b)\cdot c + (-1)^{|a|}a\cdot\Delta(b\cdot c) + 
(-1)^{(|a|+1)|b|}b\cdot \Delta(a\cdot c).
\end{multline*}
In other words, $\Delta$ is a derivation of order $2$.

Then, if we define the bilinear operation $\{,\}\colon A\otimes A\to A[-n]$ by
\[
\{a,b\}= (-1)^{|a|} \left(
\Delta(a\cdot b) - \Delta(a)\cdot b - (-1)^{|a|}a\cdot\Delta(b)
\right),
\]
we have that the quadruple $(A,\cdot,\{,\},\Delta)$ is an 
$n$-Batalin--Vilkovisky algebra. Vice versa, one easily checks that the 
product and the generator of an $n$-Batalin--Vilkovisky algebra satisfy the 
above ``seven term'' relation.

\section{Examples}

\subsection{Schouten--Nijenhuis bracket}
\label{ssec:sn}

Suppose $\mathfrak g$ is a graded Lie algebra of degree zero. Then 
$A = S(\mathfrak g [n])$ is an $n$-Poisson algebra with its natural 
multiplication (the one induced from the tensor algebra $T(A)$) and 
a degree $-n$ bracket, often called the 
Schouten--Nijenhuis bracket,  defined as follows \cite{Kosz, Kra88}
:  the bracket on $S^1(\mathfrak g [n])= \mathfrak g
[n]$ is defined as  the suspension of the bracket on
$\mathfrak g$, while on 
$S^k(\mathfrak g [n])$, for $k>1$, the bracket is defined inductively by forcing 
the Leibniz rule
\[
\{a,b\cdot c\}=\{a,b\}\cdot c + (-1)^{|b|(|a|+n)} b\cdot \{a,c\}. 
\]
Moreover, when $n$ is odd, there exists a generator defined as
\[
\Delta(a_1\cdot a_2 \cdots a_k) = \sum_{i<j} (-1)^\epsilon 
\{a_i,a_j\} \cdot a_1\cdots \widehat{a_i}\cdots \widehat{a_j}\cdots a_k
\]
where $a_1,\ldots,a_k\in\mathfrak g$ and $\epsilon=|a_i| + (|a_i|+1)
(|a_1|+\cdots |a_{i-1}|+i-1) + (|a_j|+1) (|a_1|+\cdots \widehat{|a_{i}|}+
\cdots 
|a_{j-1}|+j-2)$. An easy check shows that 
$\Delta^2=0$, thus $S(\mathfrak g [n])$ is an $n$-Batalin--Vilkovisky algebra
 for
every odd $n\in\mathbb N$. For $n=1$ the $\Delta$-cohomology on 
$\bigwedge\mathfrak g$ is the usual Cartan--Chevalley--Eilenberg cohomology.

In particular, one can consider the Lie algebra $\mathfrak g = \mathrm{Der}(B) 
= \oplus_{j\in\Z} \mathrm{Der}^j(B)$ of derivations of a graded commutative 
algebra $B$. More explicitly, $\mathrm{Der}^j(B)$ consists of linear maps 
$\phi\colon B\to B$ of degree $j$ such that $\phi(ab) = \phi(a)b+ (-1)^{j\,
 |a|} 
a\phi(b)$ and the bracket is $\{\phi,\psi\}=\phi\circ\psi 
-(-1)^{|\phi||\psi|}\psi\circ\phi$. The space of multiderivations 
$S(\mathrm{Der}(B)[1])$, endowed with the 
Schouten--Nijenhuis bracket, is a Gerstenhaber algebra.

We can further specialize to the case when $B$ is the algebra 
$\mathcal C^\infty(M)$ of smooth functions on a smooth manifold $M$; 
then $\mathfrak X(M)=\mathrm{Der}(\mathcal C^\infty(M))$ is the space
of vector fields on $M$ and 
$\mathcal V(M) = S(\mathfrak X(M))[1])$ 
is the space of multivector fields on $M$. It is a classical result 
by Koszul \cite{Kosz} that there is a bijective correspondence 
between generators for $\mathcal V(M)$ and connections on the highest
exterior power $\bigwedge^{{\rm dim}M}TM$ of the tangent bundle of $M$. 
Moreover, flat connections correspond to generators which square to zero.

\subsection{Lie algebroids}

A Lie algebroid $E$ over a smooth manifold $M$ 
is a vector bundle $E$ over $M$ together with
a Lie algebra structure (over $\mathbb R$) on the space $\Gamma(E)$ of 
smooth sections of $E$, 
and a bundle map $\rho\colon E\to TM$, called the anchor, extended to a
map $\rho_*$ between 
sections of these bundles, such that
\begin{align*}
\{X, fY \} &= f\{X, Y \} + (\rho_*(X)f)Y
\end{align*}
for any smooth sections $X$ and $Y$ of $E$ and any smooth function $f$ on $M$.
 In
particular, the anchor map on sections $\rho_*^{}\colon
\Gamma(E)\to \mathfrak X(M)$ is a morphism of Lie 
algebras , namely $\rho_*^{}(\{X, Y \}) =
\{\rho_*^{}(X), \rho_*^{}(Y)\}$.
 
The link between Lie algebroids and Gerstenhaber algebras is 
given by the following Proposition \cite{KS95b,Xu}:

\begin{prop}
Given a vector bundle $E$ over $M$, there exists a one-to-one 
correspondence between Gerstenhaber algebra structures on 
$A=\Gamma(\bigwedge(E))$ and Lie algebroid structures on $E$.
\end{prop}
 
The key of the Proposition is that one can extend the Lie algebroid 
bracket to a unique graded antisymmetric bracket on 
$\Gamma(\bigwedge(E))$ such that $\{X,f\}=\rho(X)f$ for $X\in 
\Gamma(\bigwedge^1(E))$ 
and $f\in\Gamma(\bigwedge^0(E))$, and that for $Q\in\Gamma(\bigwedge^{q+1}(E))$, 
$\{Q,\cdot\}$ is a derivation of $\Gamma(\bigwedge(E))$ of degree $q$.
 
\begin{exmp}\label{exmp:g}
A finite dimensional Lie algebra $\mathfrak g$ can be seen as a Lie
algebroid over a trivial base manifold. The corresponding Gerstenhaber algebra 
is the one of subsection~\ref{ssec:sn}.
\end{exmp} 
 
\begin{exmp}
\label{exmp:tan}
The tangent bundle $TM$ of a smooth manifold $M$ is a Lie algebroid with anchor
map given by the identity and algebroid Lie bracket given by the usual Lie
bracket on vector fields. In this case we recover the Gerstenhaber algebra of
multivector fields on $M$ described in subsection~\ref{ssec:sn}. 
\end{exmp} 

\begin{exmp}
\label{exmp:cotan}
If $M$ is a Poisson manifold with Poisson bivector field $\alpha$, then the 
cotangent bundle $T^*M$ inherits a natural Lie algebroid structure where 
the anchor map $\alpha^\#\colon T_p^*M\to T_pM$ at the point $p\in M$ is 
given by $\alpha^\#(\xi)(\eta)=\alpha(\xi,\eta)$, with $\xi,\eta\in T_p^*M$,
and the Lie bracket of the 1-forms $\omega_1$ and $\omega_2$ is given by
\[
\{\omega_1,\omega_2\} = \LL_{\alpha^\#(\omega_1)}\omega_2 -
\LL_{\alpha^\#(\omega_2)}\omega_1 - \mathrm{d}\alpha(\omega_1,\omega_2)
\]
The associated Gerstenhaber algebra is the de~Rham algebra of 
differential forms endowed with the bracket defined by Koszul 
in \cite{Kosz}. As shown in \cite{KS95b}, $\Gamma(\bigwedge(T^*M))$ 
is indeed a BV algebra with an exact generator $\Delta = 
[\mathrm{d},\iota_\alpha]$ given by the commutator of the contraction 
$\iota_\alpha$ with the Poisson bivector $\alpha$ and the de~Rham 
differential $\mathrm{d}$. Similar results hold if $M$ is a 
Jacobi manifold. 
\end{exmp} 

It is natural to ask what additional structure on a Lie algebroid $E$ makes 
the Gerstenhaber algebra $\Gamma(\bigwedge(E))$ into a BV algebra. 
The answer is given by following result, which is proved in \cite{Xu}
\begin{prop}
Given a Lie algebroid $E$, there is a one-to-one correspondence
between generators for the Gerstenhaber algebra $\Gamma(\bigwedge(E))$ 
and $E$-connections on $\bigwedge^{\mathrm{rk}E}E$ (where $\mathrm{rk}E$
denotes the rank of the vector bundle $E$). Exact generators 
correspond to flat $E$-connections, and in particular, since flat 
$E$-connections always exist, $\Gamma(\bigwedge(E))$ is always a BV algebra.
\end{prop}

For the appearance of Gerstenhaber algebras in the theory of Lie bialgebras and
Lie bialgebroids, see
\cite{L-R90}.

\subsection{Lie algebroid cohomology}\label{s:Lac}
A Lie algebroid structure on $E\to M$ defines a differential $\delta$ on 
$\Gamma(\bigwedge E^*)$ by
\[
\delta f:=\rho^* \mathrm{d} f,\quad f\in C^\infty(M)=\Gamma(\Lambda^0 E^*)
\]
and
\begin{multline*}
\braket{\delta\sigma}{X\wedge Y}:=
\braket{\delta\braket\sigma X}Y
-\braket{\delta\braket\sigma Y}X-\braket\sigma{\{X,Y\}},\\ 
X,\,Y\in\Gamma(E),\ \sigma\in\Gamma(E^*),
\end{multline*}
where $\rho^*\colon\Omega^1(M)\to\Gamma(E^*)$ is the transpose of 
$\rho_*\colon\Gamma(E)\to\mathfrak X(M)$ and
$\braket{\ }{\ }$ is the canonical pairing of sections 
of $E^*$ and $E$. On $\Gamma(\bigwedge^n E^*)$, with
$n\geq2$, the differential $\delta$ is defined by forcing
the Leibniz rule.

In example~\ref{exmp:g} we get the Cartan--Chevalley--Eilenberg differential
on $\bigwedge\mathfrak g^*$; in example~\ref{exmp:tan} we recover
the de~Rham differential on $\Omega^*(M)=\Gamma(\bigwedge T^*M)$,
while in example~\ref{exmp:cotan} the differential on 
$\mathcal V(M)=\Gamma(\bigwedge TM)$ is $\{\alpha,\ \}_{SN}$.

\subsection{Lie--Rinehart algebras}

The algebraic generalization of a Lie algebroid is a Lie--Rinehart algebra.
Recall that given a commutative associative algebra $B$ (over some ring $R$)
and a 
$B$-module $\mathfrak g$, then a Lie--Rinehart algebra structure 
on $(B,\mathfrak g)$ is a Lie algebra structure (over $R$) on $\mathfrak g$ 
and an action of $\mathfrak g$ on the left on $B$ by derivations, 
satisfying the following compatibility conditions:
\begin{align*}
\{\gamma, a\sigma\} & = \gamma(a)\sigma + a\{\gamma, \sigma\} \\
(a\gamma)(b) & = a(\gamma(b))  
\end{align*}
for every $a,b\in B$ and $\gamma,\sigma\in\mathfrak g$. 

The Lie--Rinehart structures on the pair $(B,\mathfrak g)$
bijectively correspond to the Gerstenhaber algebra
structures on the exterior algebra $\bigwedge_B(\mathfrak g)$ 
of $\mathfrak g$ in the category of $B$-modules. When $\mathfrak g$ is
of finite rank over $B$, generators  for these structures are in turn 
in bijective
correspondence  with $(B,\mathfrak g)$-connections on
$\bigwedge_B^{\mathrm{rk}_B^{}{\mathfrak g}}{\mathfrak g}$, 
and flat connections 
correspond to exact generators. 
For additional discussions,  
see \cite{G-S92, Hueb}.

Lie algebroids are Lie--Rinehart algebras in the smooth setting. Namely,
if $E\to M$ is a Lie algebroid, then the pair 
$(\mathcal C^\infty(M), \Gamma(E))$ is a Lie--Rinehart algebra (with
action induced by the anchor and the given Lie bracket).

\subsection{Lie--Rinehart cohomology}\label{s:LRc}
Lie algebroid cohomology may be generalized to every Lie--Rinehart algebra
$(B,\mathfrak g)$. Namely, on the complex is $\mathrm{Alt}_B(\mathfrak g,B)$
of alternating 
multilinear functions on $\mathfrak g$ with values in $B$, one can define
a differential $\delta$ by the rules
\[
\langle\delta a,\gamma\rangle =\gamma(a),\qquad a\in
B=\mathrm{Alt}_B^0(\mathfrak g,B),
\ \gamma\in\mathfrak g,
\]
\begin{multline*}
\langle \delta a,\gamma\wedge\sigma\rangle =
\langle\delta\langle a,\gamma\rangle,\sigma\rangle -
\langle\delta\langle
a,\sigma\rangle,\gamma\rangle-\langle
a,\{\gamma,\sigma\}\rangle ,
\quad
\gamma,\,\sigma\in\mathfrak g,\ 
a\in \mathrm{Alt}_B^1(\mathfrak g,B),
\end{multline*}
and forcing the Leibniz rule on elements of
$\mathrm{Alt}_B^n(\mathfrak g,B)$, $n\geq 2$.

\subsection{Hochschild cohomology}

Let $A$ be an associative algebra 
with product $\mu$, and consider the Hochschild cochain 
complex $\mathsf{Hoch}(A)=\prod_{n\ge 0}\mathrm{Hom}
(A^{\otimes n},A)[-n+1]$. There are two basic 
operations between two elements 
$f\in\mathrm{Hom}(A^{\otimes k},A)[-k+1]$ and 
$g\in\mathrm{Hom}(A^{\otimes l},A)[-l+1]$, namely a degree zero 
product
\[
(f\cup g)(a_1\otimes\cdots\otimes a_{k+l}) = (-1)^{kl} 
f(a_1\otimes\cdots\otimes a_l)\,
g(a_{l+1}\otimes\cdots\otimes a_{k+l})
\]
and a degree -1 bracket $\{f,g\}=f\circ g -(-1)^{(k-1)(l-1)}g\circ f$, 
where
\begin{multline*}
(f\circ g)(a_1\otimes\cdots\otimes a_{k+l-1}) =\\=
\sum_{i=1}^{k-1} (-1)^{i(l-1)} f(a_1\otimes\cdots\otimes
a_i\otimes  g(a_{i+1}\otimes\cdots\otimes a_{i_l})
\otimes \cdots\otimes a_{k+l-1}).
\end{multline*}
It is well known from \cite{Gerst} that the cohomology 
$\mathsf{HHoch}(A)$ of the Hochschild complex with 
respect to the differential
${\mathrm d}^{}_{\mathsf{Hoch}}=\{\mu,\cdot\}$ has the structure of a 
Gerstenhaber 
algebra. More generally, there is a Gerstenhaber algebra 
structure on Hochschild cohomology of differential graded
associative algebras \cite{Loday} and of  $A_\infty$
algebras. Moreover, if $A$ is endowed with a symmetric, invariant and non-degenerate  inner
product, then $\mathsf{HHoch}(A)$ is also a BV algebra \cite{Tradler}. 

\subsection{Graded symplectic manifolds}

The construction of Example~\ref{exmp:smoothnocoor} can be extended to 
graded symplectic manifolds; see \cite{AKSZ,Getz,Schw93}. 
Recall that a symplectic structure of degree $n$ on a graded
manifold  
$N$ is a closed nondegenerate two-form $\omega$ such that 
$L_E\omega=n\omega$ where $L_E$ is the Lie derivative with respect to 
the Euler field of $N$ (see \cite{Roy02b} for details). 
Let us denote by $X_h$ the vector field associated to the 
function $h\in \mathcal C^\infty(N)$ 
by the formula $\iota_{X_h}\omega=\mathrm{d}h$. Then the bracket 
\[
\{f,g\} = \iota_{X_f}\iota_{X_g} \omega
\]
gives $\mathcal C^\infty(N)$ the structure of a graded Poisson algebra of 
a degree $n$.

If the symplectic form has odd degree and the graded manifold has a 
volume form, then it is possible to construct an exact generator defined 
by $$\Delta(f)=\frac12 \mathrm{div}(X_f)$$ where $\mathrm{div}$ is the 
divergence operator associated to the given volume form 
\cite{Getz,KS-M02}.

An explicit characterization of graded symplectic manifolds 
has been given in \cite{Roy02b}.
In particular it is proved there that every symplectic form of degree $n$ with
$n\ge 1$ is necessarily exact. More precisely, one has 
$\omega=\mathrm d(\iota_E\omega/n)$.

\subsection{Shifted cotangent bundle}
\label{ssec:cotan}

The main examples of graded symplectic manifolds are given by shifted 
cotangent bundles. If $N$ is a graded manifold then the shifted cotangent 
bundle $T^*[n]N$ is the graded manifold obtained by shifting by $n$ the 
degrees of the fibers of the cotangent bundle of $N$. This graded manifold 
possesses a non degenerate closed two-form of degree $n$, 
which can be expressed in local coordinates as
\[
\omega = \sum_i \mathrm{d}x^i \wedge \mathrm{d}x^\dagger_i
\]
where $\{x^i\}$ are local coordinates on $N$ and $\{x^\dagger_i\}$ 
are coordinate functions on the fibers of $T^*[n]N$.
In local coordinates, the bracket between two homogeneous functions 
$f$ and $g$ is given by
\[
\{f,g\} = -(-1)^{|x^\dagger_i|\,|f|}
 \frac{\partial f}{\partial x^\dagger_i}\ 
\frac{\partial g}{\partial x^i} 
-(-1)^{(|f|+n)(|g|+n) + |x^\dagger_i|\,|g|}
\frac{\partial g}{\partial x^\dagger_i}\ 
\frac{\partial f}{\partial x^i} 
\]
If in addition the graded manifold $N$ is orientable, then $T^*[n]N$ 
has a volume form too; when $n$ is odd, the exact generator $\Delta(f) =
\frac12 
\mathrm{div} X_{f}$, already considered in \cite{Kosz},
 is written in local coordinates as 
\[
\Delta = \frac{\partial}{\partial x^\dagger_i}
\frac{\partial}{\partial x^i}.
\]

In the case $n=-1$, we have a natural identification between functions
on $T^*[-1]N$ and multivector fields $\mathcal V(N)$ on $N$
and we recover again the Gerstenhaber algebra of subsection~\ref{ssec:sn}.
Moreover it is easy to see that, under the above identification, 
$\Delta$ applied to a vector field of $N$ is the usual divergence 
operator.

\subsection{Examples from algebraic topology}

For any $n>1$, the homology of the $n$-fold loop space 
$\Omega^n(M)$ of a topological space $M$ has the structure 
of an $(n-1)$-Poisson algebra \cite{May}. In particular the 
homology of the double loop space $\Omega^2(M)$
is a Gerstenhaber algebra, and has an exact generator defined using 
the natural circle action on this space \cite{Getz}. The  homology 
of the free loop space $\mathcal L(M)$ of a closed oriented manifold $M$ is 
also a
BV algebra when endowed with the ``Chas--Sullivan intersection product'' 
and with a generator defined again using the natural 
circle action on the free loop space. Recently it has been 
shown that the homology of the space of framed embeddings of $\mathbb R^n$ 
into $\mathbb R^{n+k}$ is an $n$-Poisson algebra.

\subsection{Shifted Poisson structures on derived stacks}
In recent years, a vast generalization of the notion of Poisson structure on a smooth manifold has been introduced in the context of derived algebraic geometry by Calaque, Pantev, T\"oen, Vaqui\'e and Vezzosi who considered shifted Poisson structures on derived Artin stacks \cite{cptvv}. Basic examples are the 1-shifted and 2-shifted Poisson structures on
$BG$ given by elements in $(\wedge^3\mathfrak{g})^{\mathfrak{g}}$ and in $(S^2\mathfrak{g})^\mathfrak{g}$, respectively, where $G$ is a reductive algebraic group and $\mathfrak{g}$ is its Lie algebra.
The corresponding notion for smooth stacks presented by Lie groupoids has then been investigated by Bonechi, Ciccoli, Laurent-Gengoux and Xu in \cite{bclgx}, and by Ginot, Ortiz and Stefani in \cite{gos}.
\par
One starts by considering, for any $m\in \mathbb{Z}$ the differential graded commutative algebra
\[
{\sf{MultiVect}}(X,m)=\bigoplus_{p\geq 0} \Gamma(X,S^p(\mathbb{T}_X[-m])
\]
of $m$-shifted multivector fields on $X$. Here $\mathbb{T}_X$ denotes the tangent complex of the stack $X$ and $\Gamma$ is the derived global sections (i.e., hypercohomology) functor. The algebra ${\sf{MultiVect}}(X,m)$ is bigraded, with one grading given by the cohomological degree and the other, called the weight given by the exponent $p$. The main result is then that ${\sf{MultiVect}}(X,m)$ carries a natural $\mathbb{P}_{m+1}$-commutative differential graded algebra structure, where $\mathbb{P}_{m+1}$ is the dg-operad whose algebras are Poisson cdga's with a bracket of degree $-m$. 
Equivalently, (commutative) $\mathbb{P}_{m+1}$-algebras can be inductively defined as (commutative) associative algebras in $\mathbb{P}_{m}$-algebras \cite{Safronov}. 
Moreover, the Poisson bracket on $\mathbb{P}_{m+1}$ has weight $-1$. Therefore, ${\sf{MultiVect}}(X,m)[m]$ is a dg-Lie algebra with a Lie bracket of weight $-1$, and one can define the space ${\sf{Poiss}}(X,n)$ of $n$-shifted Poisson structures on $X$ as
\[
{\sf{Poiss}}(X,n)={\sf{Map}}_{{\sf{dgLie}}^{\sf{gr}}}(\mathbb{K}(2)[-1],{\sf{MultiVect}}(X,n+1)[n+1]),
\]
where ${\sf{Map}}_{{\sf{dgLie}}^{\sf{gr}}}$ is the Hom-space in the $(\infty,1)$-category of graded dg-Lie algebras over the (characteristic zero) field of definition $\mathbb{K}$ of the stack $X$, and $\mathbb{K}(2)[-1]$ is the graded dg-Lie algebra consisting of $\mathbb{K}$ in pure cohomological degree 1, pure weight 2, and trivial bracket. When $X$ is a smooth underived scheme, ${\sf{Poiss}}(X,0)$ is the usual set of Poisson bivectors on $X$. Moreover, one has a natural notion of nondegenerate shifted Poisson strutures and an equivalence ${\sf{Poiss}}(X,n)^{\sf{nd}}\simeq {\sf{Symp}}(X,n)$ between the space of nondegenerate $n$-shifted Poisson structures on $X$ and the space of $n$-shifted symplectic structures on $X$ as defined in \cite{ptvv}

\section{Applications}

\subsection{BRST quantization in the Hamiltonian formalism}

The BRST procedure, after Becchi--Rouet--Stora and Tyutin \cite{BRST,BRST2,Tyu75}, is a method for quantizing classical mechanical 
systems or classical field theories in the presence 
of symmetries. We describe here its classical, Hamiltonian counterpart 
following \cite{Sta97}.
The starting point is a symplectic manifold $M$ (the ``phase 
space''), a function $H$ (the ``Hamiltonian'' of the system) 
governing the evolution of the system, and the 
``constraints''
given by several functions $g_i$ which
Poisson commute with $H$ and among each other up to a 
$\mathcal C^\infty(M)$-linear combination of the $g_i$'s.

Then 
the dynamics is constrained on the locus $V$ of common zeroes 
of the $g_i$'s. When $V$ is a submanifold, the  
$g_i$'s are a set of generators for the
ideal $I$ of functions vanishing on $V$.
Observe that $I$ is
closed under the Poisson bracket. Functions in $I$ are 
called ``first class constraints.'' The Hamiltonian vector fields
of first class constraints, which are by construction tangential to
$V$, are the ``symmetries'' of the system.

When $V$ is smooth, then 
it is a coisotropic submanifold of $M$ and the 
Hamiltonian vector fields determined by the constraints 
give a foliation $\mathcal F$ of $V$. 
In the nicest case $V$ is a
principal bundle with $\mathcal F$ its vertical foliation 
and the algebra of functions $\mathcal C^\infty(V/\mathcal F)$ 
on the ``reduced phase space''  $V/\mathcal F$ is identified
with  the $I$-invariant subalgebra of $\mathcal C^\infty(M)/I$. 

{}From a
physical point of view, the points of
$V/\mathcal F$ are  the interesting states at a classical level, and
a quantization  of this system means a quantization of $\mathcal
C^\infty(V/\mathcal F)$.  The BRST procedure gives a method of
quantizing 
$\mathcal C^\infty(V/\mathcal F)$ starting from the (known) 
quantization of $\mathcal C^\infty(M)$. Notice that these 
notions immediately generalize to graded symplectic 
manifolds.

{}From an algebraic point of view, one starts with a graded 
Poisson algebra $P$ and a multiplicative ideal $I$ which is closed 
under the Poisson bracket. The algebra of functions on the 
``reduced phase space'' is replaced by $(P/I)^I$, 
the $I$-invariant subalgebra of $P/I$. This 
subalgebra inherits a Poisson bracket even if $P/I$ does not. 
Moreover the pair $(B,\mathfrak g)=(P/I,I/I^2)$ inherits a graded 
Lie--Rinehart structure.
The ``Rinehart
complex'' $\mathrm{Alt}_{P/I}(I/I^2,P/I)$ of alternating 
multilinear functions on $I/I^2$ with values in $P/I$, endowed
with the differential described in~\ref{s:LRc}, plays the role of the 
de~Rham complex of vertical forms on $V$ with respect to the 
foliation $\mathcal F$ determined by the constraints. 

In case $V$ is a smooth submanifold, we also have the following geometric
interpretation: Let $N^*V$ denote the conormal bundle of $V$ (i.e., the
annihilator of $TV$ in $T_VP$). This is a Lie subalgebroid of $T^*P$ 
if and only
if $V$ is coisotropic. Since we may identify $I/I^2$ with sections of
$N^*C$ (by the de~Rham differential), $(P/I,I/I^2)$ is the corresponding
Lie--Rinehart pair. The Rinehart complex is then the corresponding Lie 
algebroid complex $\Gamma(\bigwedge (N^*V)^*)$ with differential described
in~\ref{s:Lac}. The image of the anchor map $N^*V\to TV$ is the distribution
determining $\mathcal F$, so by duality we get an injective chain map
form the vertical de~Rham complex to the Rinehart complex.

The main point of the BRST procedure is to define a chain complex 
$C^\bullet = \bigwedge (\Psi^*\oplus\Psi) 
\otimes P$,
where $\Psi$ is a graded vector space, with a coboundary 
operator $D$ (the ``BRST operator''), and a quasi-isomorphism (i.e., a chain map
that induces an isomorphism in cohomology)
\[
\pi\colon (C^\bullet,D)\to  \left(\mathrm{Alt}_{P/I}(I^2/I,P/I),
\mathrm{d}\right).
\]
This means in particular that the zeroth cohomology $H_D^0(C)$ gives 
the algebra $(P/I)^I$ of functions on the ``reduced phase space''.
Observe that there is a natural symmetric inner product on $\Psi^*
\oplus\Psi$ given by the evaluation of $\Psi^*$ 
on $\Psi$. This inner product, as an element of $S^2(\Psi
\oplus\Psi^*)\simeq S^2(\Psi)\oplus (\Psi
\otimes\Psi^*)\oplus S^2(\Psi^*)$, is concentrated in the component $\Psi
\otimes\Psi^*$, and so it defines an element in $\bigwedge^2(\Psi[1]
\oplus\Psi^*[-1])\simeq S^2(\Psi)[2]
\oplus(\Psi\otimes\Psi^*)
\oplus S^2(\Psi^*)[-2]$, i.e., a degree zero bivector field on $\Psi[1]
\oplus\Psi^*[-1]$. It is easy to see that this bivector field induces a degree
zero Poisson structure on $S(\Psi^*[-1]\oplus\Psi[1])$. {}From another viewpoint
this is the Poisson structure corresponding to the canonical symplectic structure on
$T^*\Psi[1]$. Finally, we have that 
$S(\Psi^*[-1]\oplus\Psi[1]) 
\otimes P$ is a degree zero Poisson algebra. Note that the superalgebra
underlying the graded algebra $S(\Psi^*[-1]\oplus\Psi[1])\otimes P$ is
canonically isomorphic to the complex $C^\bullet = \bigwedge
 (\Psi^*\oplus\Psi) 
\otimes P$.
When $P=\mathcal C^\infty(M)$, we can think of 
$S (\Psi^*[-1]\oplus \Psi[1]) \otimes \mathcal C^\infty(M)$ 
as the 
algebra of functions on the 
graded
symplectic 
manifold
$N=(\Psi[1]\oplus \Psi^*[-1])\times M$ (the ``extended phase 
space'').
In physical language, coordinate functions on 
$\Psi[1]$ are called ``ghost fields'' while coordinate functions 
on $\Psi^*[-1]$ are called ``ghost momenta'' or, by some 
authors, ``antighost fields''
(not to
be confused with the antighosts of the Lagrangian functional-integral
approach to quantization).

Suppose now that there exists an element $\Theta\in S
(\Psi^*[-1]\oplus\Psi[1]) \otimes P$ such that $\{\Theta,\cdot\} = D$, 
that one can extend the ``known'' quantization of $P$ to a 
quantization of $S (\Psi^*[-1]\oplus\Psi[1]) \otimes P$ as 
operators on some (graded) Hilbert space $\mathcal T$ and that the 
operator $Q$ which quantizes $\Theta$ has square zero. Then one 
can consider the ``true space of physical states'' 
$H^0_Q(\mathcal T)$ on which the $\mathrm{ad}_Q$-cohomology of operators will act.
This provides one with a quantization of $(P/I)^I$.

For further details on this procedure, and in particular for 
the construction of $D$, we refer to \cite{H-T,K-S,Sta97} 
and references therein. Observe that some authors refer to this method
as BVF [Batalin--Vilkovisky--Fradkin] and reserve the name
BRST to the case when the $g_i$'s are the components of an equivariant
moment map.

For a generalization to graded 
manifolds different from $(\Psi[1]\oplus 
\Psi^*[-1])\times M$ we refer to \cite{Roy02b}. There it is proved
that the element $\Theta$ exists if the graded symplectic form has
degree different from $-1$.

\subsection{BV quantization in the Lagrangian formalism}
\label{ssec:bvq}

The BV formalism is a procedure for the quantization of
physical systems with symmetries in the Lagrangian
formalism, see \cite{B-V,H-T}.  As a first step, the ``configuration space'' $M$
of the system is  augmented by the introduction of
``ghosts''. If $G$ is the group of  symmetries, this means
that one has to consider the graded manifold 
$W=\mathfrak g[1]\times M$. The second step is to double this space by 
introducing ``antifields for fields and ghosts'', namely one has to 
consider the ``extended configuration space'' $T^*[-1]W$, whose space of 
functions is a BV algebra by subsection~\ref{ssec:cotan}. 
The algebra of ``observables'' of this physical system is by definition 
the cohomology $H^*_\Delta (\mathcal C^\infty(T^*[-1]W))$ with respect to the 
exact generator $\Delta$.

\section{Related topics}

\subsection{AKSZ}
The graded manifold $T^*[-1]W$ considered in \ref{ssec:bvq} is a
particular example of a $QP$-manifold, i.e., of a graded
manifold $M$ endowed with an integrable (i.e., selfcommuting) vector field $Q$ of degree $1$
and a graded
$Q$-invariant symplectic structure $P$. In quantization
of classical mechanical theories, the graded symplectic manifold of
interest is the space of fields and antifields with symplectic form
of degree $1$, while $Q$ is the Hamiltonian
vector field defined by the action functional $S$; the
integrability of $Q$ is equivalent to the classical master
equation $\{S,S\}=0$ for the action functional. Quantization of
the theory is then reduced to the computation of the functional
integral $\int_{\mathcal L}{\rm exp}(iS/\hbar)$, where $\mathcal L$ is a
Lagrangian submanifold of $M$. This functional integral actually
depends only on the homology class of the Lagrangian. Locally, a
$QP$-manifold is a shifted cotangent bundle $T^*[-1]N$ and a
Lagrangian submanifold is the graph of an exact $1$-form. In the
notations of \ref{ssec:cotan}, a Lagrangian submanifold ${\mathcal
L}$ is therefore locally defined by equations
$x_i^\dagger=\partial\Phi/\partial x^i$, and the function $\Phi$
is called a fixing fermion. The action functional of interest is
then the gauge-fixed action $S\bigr\vert_{\mathcal
L}=S(x^i,\partial\Phi/\partial x^i)$.
\par
The language of $QP$-manifolds has powerful applications to
sigma-models: if $\Sigma$ is a
finite-dimensional graded manifold equipped with a volume element,
and $M$ is a $QP$ manifold, then the graded manifold
$C^\infty(\Sigma,M)$ of smooth maps from
$\Sigma$ to $M$ has a natural structure of $QP$-manifold which
describes some field theory if one arranges for the symplectic structure to
be of degree $1$. As an
illustrative example, if $\Sigma=T[1]X$, for a compact oriented
3-dimensional smooth manifold $X$, and $M={\mathfrak g}[1]$,
where $\mathfrak g$ is the Lie algebra of a compact Lie group,
the $QP$-manifold $C^\infty(\Sigma,M)$ is relevant to
Chern-Simons theory on $X$. Similarly, if $\Sigma=T[1]X$, for a
compact oriented 2-dimensional smooth manifold $X$ and $M=T[1]N$
is the shifted tangent bundle of a symplectic manifold, then the
$QP$-structure on $C^\infty(\Sigma,M)$ is related to the A-model
with target $N$; if the symplectic manifold $N$ is of the form
$N=T^*K$ for a complex manifold $K$, then one can endow
$C^\infty(\Sigma,M)$ with a complex $QP$-manifold structure,
which is related to the B-model with target $K$; this shows
that, in some sense, the B-model can be obtained from the A-model
by ``analytic continuation'' \cite{AKSZ}. 
If $\Sigma=T[1]X$, for a
compact oriented 2-dimensional smooth manifold $X$ and $M=T^*[-1]N$
with canonical symplectic structure, then the
$QP$-structure on $C^\infty(\Sigma,M)$ is related to the Poisson sigma model
($QP$-structures on $T^*[-1]N$ with canonical symplectic structure are in
one-to-one correspondence with Poisson structures on $N$).
The study of $QP$
manifolds is sometimes referred to as ``the AKSZ formalism''. 
In \cite{Roy02b} $QP$-manifolds with symplectic structure of degree $2$ are
studied and shown to be in one-to-one correspondence with 
Courant algebroids.

\subsection{Graded Poisson algebras from cohomology of $P_\infty$} 
The Poisson bracket on a Poisson manifold can be derived from the Poisson
bivector field $\alpha$ using the Schouten--Nijenhuis bracket
as follows \cite{KS96}:
\[
\{f,g\}=\{\{\alpha,f\}_{SN},g\}_{SN}.
\]
This may be generalized \cite{Vor02}
 to the case of a graded manifold $M$ endowed
with a multivector field $\alpha$ of total degree $2$ (i.e., 
$\alpha=\sum_{i=0}^\infty \alpha_i$, where $\alpha_i$ is an $i$-vector field
of degree $2-i$) satisfying the equation $\{\alpha,\alpha\}_{SN}=0$.
One then has the derived multibrackets
\begin{gather*}
\lambda_i\colon A^{\otimes i}\to A,\\
\lambda_i(a_1,\dots,a_i):=
\{\{\dots\{\{\alpha_i,a_1\}_{SN},a_2\}_{SN}\dots\}_{SN},a_i\}_{SN},
\end{gather*}
with $A=\mathcal C^\infty(M)$. Observe that
$\lambda_i$ is a multiderivation of degree $2-i$. 
The operations $\lambda_i$ define the structure
of an $L_\infty$-algebra on $A$. 
Such a 
structure is called a
$P_\infty$-algebra ($P$ for Poisson) since
the $\lambda_i$'s are multiderivations. If $\lambda_0=\alpha_0$ vanishes,
then $\lambda_1$ is a differential, and the $\lambda_1$-cohomology
inherits a graded Poisson algebra structure. This structure 
can be
used 
to describe the deformation quantization 
of coisotropic submanifolds and 
to describe their deformation 
theory.

\begin{Ack}
We thank Johannes Huebschmann, Yvette Kosmann-Schwarz\-bach,
Dmitry Roytenberg, and Marco Zambon
for very useful comments and suggestions.
\end{Ack}


\def\cprime{$'$} \def\cprime{$'$}

\end{document}